\documentclass[12pt]{article}

\title{Two-stage decision making with set relations\footnote{This article is a preprint version; it has been accepted for publication in Decisions for Economics and Finance.}}

\author{
Andreas H. Hamel\footnote{Free University of Bozen, Faculty of Economics and Management, \href{mailto:andreas.hamel@unibz.it}{andreas.hamel@unibz.it}}
}
\date{}

\usepackage{array} 
\usepackage{verbatim} 
\usepackage{float, graphicx, caption, subcaption}

\usepackage{amsmath, amssymb, enumerate, color}

\usepackage{algpseudocode,algorithm,algorithmicx}

\newtheorem{theorem}{Theorem}
\newtheorem{corollary}[theorem]{Corollary}
\newtheorem{remark}[theorem]{Remark}
\newtheorem{lemma}[theorem]{Lemma}
\newtheorem{definition}[theorem]{Definition}
\newtheorem{proposition}[theorem]{Proposition}
\newtheorem{example}[theorem]{Example}

\numberwithin{equation}{section}  
\numberwithin{figure}{section}    
\numberwithin{table}{section}     
\numberwithin{theorem}{section}

\newcommand{\of}[1]{\ensuremath{\left( #1 \right)}}

\newcommand{\cb}[1]{\ensuremath{ \left\{ #1 \right\} }}
\newcommand{\sqb}[1]{\ensuremath{ \left[ #1 \right] }}

\newcommand{\bs}{\backslash}

\newcommand{\pend}{ \hfill $\square$ \medskip}

\newcommand{\vn}{\varnothing}

\renewcommand{\P}{\ensuremath{\mathcal{P}}}

\newcommand{\A}{\ensuremath{\mathcal{A}}}
\newcommand{\B}{\ensuremath{\mathcal{B}}}
\newcommand{\F}{\ensuremath{\mathcal{F}}}

\newcommand{\MIN}{{\rm Min\,}}
\newcommand{\MAX}{{\rm Max\,}}

\newcommand{\R}{\mathrm{I\negthinspace R}}

\newcommand{\dom}{{\rm dom \,}}

\newcommand{\co}{{\rm co \,}}

\newcommand{\triup}{{\rm \vartriangle}}
\newcommand{\trido}{{\rm \triangledown}}

\usepackage[
colorlinks=true, linkcolor=blue, citecolor=blue, urlcolor=blue, filecolor=blue
auskommentieren und unten pdfborder aktivieren]{hyperref}

\definecolor{color0}{gray}{.50}
\definecolor{color1}{rgb}{0,.2,.8}
\definecolor{color2}{rgb}{1,.2,0}
\definecolor{color3}{rgb}{.8,.5,1}
\newcommand{\f}{\color{color1}}

\begin{document}

\maketitle

\begin{abstract}
Two-stage decision problems with incomplete preferences associated to a set-valued function are studied. Depending on the preferences of the decision makers collaborative and contradictive decision problems are distinguished. Both types are embedded into set optimization problems with objective functions mapping into complete lattices of sets. Structural properties of these lattices are given and decision making procedures are outlined. The case of partially aligned preferences is discussed. Applications to a game theoretic situation including new equilibrium concepts are given.
\end{abstract}

{\bf Keywords.} decision problem, incomplete preference, set relations, complete lattice

\medskip
{\bf MSC.} 90B50, 91B06

\section{Introduction}

Kreps in \cite{Kreps79Economet} discussed a decision making procedure which has a two-stage character where in the first stage a 'menu,' i.e., a set of alternatives, has to be chosen and in the second stage, only later, a 'dish,' i.e., an element of the previously selected set. An important feature is that the decision maker might have a 'preference for flexibility' which means that the first stage decision should leave as many options as possible to the second stage decision. This is also called the Cautios Deferral rule in the economic literature, compare \cite{PejsachowiczToussaert17JET}.

Moreover, Kreps based the first stage decision on preference relations for sets which are nowadays known as set relations and form the very basis of set optimization \cite{HamelEtAl15Incoll}, a new branch of optimization theory, with many applications, e.g., in finance \cite{HamelHeyde21Math} (see the references therein), statistics \cite{HamelKostner18JMVA}, game theory \cite{HamelLoehne18MMOR}, robust optimization \cite{EhrgottEtAl14EJOR, IdeEtAl14FPTA}, multi-criteria decision making \cite{ZitzlerThieleBader09IEEE} (with precursors dating back to 1999) and vector optimization \cite{KuroiwaTanakaTruong97JNCA} with many follow-ups which popularized set relations within the  multi-objective and vector optimization community.

The completeness of the order relation for single items was assumed in \cite{Kreps79Economet}. This assumption was later removed, e.g., in \cite{ArlegiNieto01MSS, PejsachowiczToussaert17JET}. 

In this paper, a stringent set optimization approach to the two-stage decision problem with {\f incomplete} preferences is developed. Elements of this approach can already been found in \cite{HamelLoehne24ArX, Loehne24ArX}.

Main contributions of the paper include: (1) a theoretical framework for two-stage decision making with incomplete preferences for flexibility, (2) complete lattice-valued set optimization problems and corresponding solution concepts upon which decision making procedures can be based, (3) characterization of solutions involving sets of desirable and acceptable outcomes, respectively, (4) concepts for situations with aligned and inverse preferences as well as partially aligned preferences, (5) solution concepts for games with general payoffs for which only incomplete preferences are available.

The plan of the paper is as follows. The next section introduces two versions of the basic problem, the collaborative and the contradictive case. Basics on relations and their extensions to set relations are collected in Section \ref{SecBasic} while in Section \ref{SecBasicDMs} the models are refined and some exemplary decision situations are described. Section \ref{SecSetOptDM} is the core part of the paper in which the set optimization approach to two-stage decision making is developed. Moreover, decision making procedures are outlined and discussed which might be considered as deriving managerial guidelines from the 'preference for flexibility' approach.

In the following section,  the "neither collaborative, nor contradictive" case is discussed. In Section \ref{SecGames}, an application to games with vector payoffs are considered and two equilibrium definitions are given generalizing those from \cite{HamelLoehne18MMOR}: one is a worst-case equilibrium, the other one of Nash type. Finally, a few conclusions and outlooks to further research topics are given in Section \ref{SecConclusions}.

\section{Decision problems associated to set-valued functions}

Let $X$, $Z$ be two non-empty sets and $\P(Z)$ the power set of $Z$, i.e., the set of all subsets of $Z$ including $\vn$.

Let $f \colon X \to \P(Z)$ be a (set-valued) function. The decision making problem associated with $f$ is a two-stage selection procedure.

First, select $x \in X$. Secondly (maybe later), select $z \in f(x)$.

Depending on by whom these choices are made, two problem classes are defined as follows.

\begin{itemize}

\item[(Coll)] The first and the second stage decision(s) are made by the same decision maker or by two decision makers whose interests are totally aligned. The problem is called the {\em collaborative decision making problem} associated to $f$.

\item[(Con)] The interests of the second stage decision maker run totally against the ones of the first stage decision maker. The problem is called the {\em contradictive decision making problem} associated to $f$.

\end{itemize}
Moreover, a third case will be considered in which the decision makers share some interest but are not totally aligned. To the best of my knowledge, this case and also the contradictive one have not been discussed in the literature on preferences for flexibility or set rankings.

\section{Basic notation}
\label{SecBasic}

A (binary) relation on a set $Z$ is a non-empty set $R \subseteq Z \times Z$, i.e., a collection of (ordered) pairs of elements of $Z$. The diagonal relation is $D = \cb{(z, z) \mid z \in Z}$. The composition of two relations $R, S \subseteq Z \times Z$ is the relation
\[
S \circ R = \cb{(u, w)\in Z \times Z \mid \exists v \in Z \colon (u, v) \in R, (v, w) \in S}.
\]
A relation $R$ is reflexive if $D \subseteq R$. It is transitive if $R \circ R \subseteq R$. A reflexive and transitive relation is called a preorder in which case $(u, v) \in R$ is written as $u \preceq v$. A relation $R$ is called complete if $(u, v) \in R$ or $(v, u) \in R$ or both for every $u, v \in Z$. Otherwise, it is called incomplete.

The inverse of a relation $R$ is the relation $R^{-1} = \cb{(v, u) \mid (u, v) \in R}$. A relation is symmetric if $R = R^{-1}$, and it is antisymmetric if $R \cap R^{-1} \subseteq D$. An antisymmetric preorder is called a partial order while a symmetric preorder is called an equivalence relation. The symmetric part of a relation $R$ is $S(R) = R \cap R^{-1} = \{(u, v) \in Z \times Z \mid (u, v) \in R, (v, u) \in R\}$.

If $E \subseteq Z \times Z$ is an equivalence relation, then $[z]_E = \{u \in Z \mid (u, z) \in S(R)\}$ is the equivalence class generated by $z \in Z$.

Let $\preceq$ be a partial order. The infimum of a set $V \subseteq Z$ is an element $\bar z \in Z$ satisfying (i) $\bar z \preceq z$ for all $z \in V$ and (ii) $w \preceq \bar z$ whenever $w \preceq z$  for all $z \in V$. If the infimum exists, it is unique and denoted by $\inf V$. Likewise, the supremum $\sup V$ is defined. A partially ordered set is called a complete lattice if each set has an infimum and a supremum within the set. Since this applies in particular to the whole set, a complete lattice always has a least (bottom) and a greatest (top) element which are both unique. Note that the order in a complete lattice does not need to be a complete one.

Let $\preceq$ be a preorder and $V \subseteq Z$. An element $\bar z \in V$ is a minimal element of $V$ if
\[
z \in V, z \preceq \bar z \; \Rightarrow \; \bar z \preceq z.
\]
The set of minimal elements of $V$ is denoted by $\MIN V$ and can be empty, a singleton or a non-singleton set. Likewise, the set of maximal elements $\MAX V$ is defined. If $\preceq$ is also antisymmetric, i.e., a partial order, $\bar z \preceq z$ in the above formula can be replaced by $\bar z = z$ due to antisymmetry.

A preorder $R$ on $Z$ can be extended to a preorder on $\P(Z)$ in two different ways. Indeed, two binary relations on $\P(Z)$ can be introduced by
\begin{align*}
(A, B) \in R^\triup \quad & \Leftrightarrow \quad \forall b \in B, \exists a \in A \colon (a, b) \in R \\
(A, B) \in R^\trido \quad & \Leftrightarrow \quad \forall a \in A, \exists b \in B \colon (a, b) \in R
\end{align*}
for $A, B \in \P(Z)$. These set relations have been (re-)introduced independently of each other several times in the literature. The article \cite{Kreps79Economet} (with the assumption that the base relation $R$ is complete) is an example from the economic literature, \cite[3.1 Definition]{Brink93AU} (the intersection of $R^\triup$ and $R^\trido$, see also 3.3) is one within a purely algebraic framework and \cite{KuroiwaTanakaTruong97JNCA} one with vector optimization in view. There are many more references, see \cite[Sect. 2.2]{HamelEtAl15Incoll} for an overview.

As one easily checks, $R^\triup$ and $R^\trido$ are preorders which are not partial orders in general even if $R$ is a partial order on $Z$. Moreover, the two relations $R^\triup$ and $R^\trido$ do not coincide in general, and one has
\begin{align}
\label{EqPfFMin}
A \supseteq B \quad & \Rightarrow \quad (A, B) \in R^\triup \\
\label{EqPfFMax}
A \subseteq B \quad & \Rightarrow \quad (A, B) \in R^\trido
\end{align}
while $A \supseteq B$ does not imply $(A, B) \in R^\trido$ and $A \subseteq B$ does not imply $(A, B) \in R^\triup$.

Finally, both relations extend $R$ from $Z$ to $\P(Z)$ in the following sense: for $u, v  \in Z$ one has
\[
(u , v) \in R \quad \Leftrightarrow \quad \of{\{u\}, \{v\}} \in R^\triup \quad \Leftrightarrow \quad
	\of{\{u\}, \{v\}} \in R^\trido.
\]
In the following, $A \preceq^\triup B$ is written for $(A, B) \in R^\triup$ and $A \preceq^\trido B$ for $(A, B) \in R^\trido$. 

Note that \eqref{EqPfFMin} expresses that decision makers with $\preceq^\triup$ prefer to have more options than less if minimization with respect to $R^\triup$ is the goal while \eqref{EqPfFMax} has the same interpretation with respect to maximization with respect to $R^\trido$. This feature has been called 'preference for flexibility' in \cite{Kreps79Economet}.

Examples mostly involve the $d$-dimensional space $Z = \R^d$ of vectors with $d$ real components. In such a case, an element $z \in \R^d$ is also written $z = (z_1, \ldots, z_d)$ where the $z_i$'s are the components of $z$. A convex cone $C \subseteq \R^d$ is a set satisfying $s C \subseteq C$ for all $s > 0$ and $C + C \subseteq C$. Such a cone generates a preference $\leq_C$ by $y \leq_C z$ iff $z -y \in C$ for $y, z \in \R^d$. The relation $\leq_C$ is transitive since $C$ is a convex cone and reflexive iff $0 \in C$. In this case, it also is a vector preorder, i.e., compatible with the algebraic operations on $\R^d$. A vector preorder is antisymmetric, i.e., a partial order, iff $C \cap (-C) = \{0\}$. The component-wise order on $\R^d$ is generated by the convex cone $\R^d_+$ of vectors with only non-negative components.

\section{Preferences of the decision makers}
\label{SecBasicDMs}

Let $R_1$, $R_2$ (also written $\preceq_1$, $\preceq_2$) be two preorders on $Z$. It is assumed that the first stage decision maker has preference $R_1$ for alternatives in $Z$ and the second stage decision maker has preference $R_2$.

\begin{definition}
Two decision makers DM1, DM2 are called {\em collaborative} if $R_2 = R_1$, i.e.,
\[
z  \preceq_1 z' \quad \Leftrightarrow \quad z \preceq_2 z' 
\]
for $z, z' \in Z$, i.e., the two preferences coincide. Two decision makers DM1, DM2 are called {\em contradictive} if $R_2 = (R_1)^{-1}$, i.e., for $z, z' \in Z$ one has
\[
z  \preceq_1 z' \quad \Leftrightarrow \quad z'  \preceq_2 z.
\]
\end{definition}

Instead of considering two DMs with respective preferences, one can also understand the situation as a dynamic one: DM1 is the first stage, DM2 the second stage decision maker. The collaborative case then means that the preference of the decision maker does not change after having made the first stage decision while the preference would completely turn around in the contradictive case. It is assumed throughout that the decision makers prefer $z \in Z$ over $z' \in Z$ if $z \preceq_1 z'$ and $z \preceq_2 z'$, respectively, i.e., the point of view is minimality (minimize losses).

\subsection{Aligned preferences}

In the collaborative case, it can be assumed that the two decisions are made by the same agent, denoted DM. The common preference is denoted by $\preceq$. The question now is how to select $x \in X$ and then how to select $z \in f(x)$.

Assume the first stage decision has been made, i.e., $x \in X$ is selected. Then it does not make sense for DM to pick $z \in f(x)$ such that there  is $z' \in f(x)$ with $z' \preceq z$ and $z \not\preceq z'$. This means that DM should pick $\bar z \in \MIN f(x)$, i.e., a minimal element in $f(x)$ since minimization is the goal. Of course, this makes only sense if $\vn \neq \MIN f(x)$.

Next, assume that DM has to compare $x, x' \in X$ for a first stage decision. The alternative $x'$ is certainly less appealing than $x$ if
\[
\forall z' \in f(x'), \; \exists z \in f(x) \colon z \preceq z',
\]
i.e., $f(x) \preceq^\triup f(x')$. DM would at least not be worse off by selecting $x$ instead of $x'$, maybe even better.

Therefore, DM should look for an alternative $\bar x \in X$ such that there is no $x \in X$ satisfying $f(x) \preceq^\triup f(\bar x)$, but not $f(\bar x) \preceq^\triup f(x)$. Thus, the DM should look for $\bar x \in X$ such that $f(\bar x)$ is a minimal element with respect to the set relation $\preceq^\triup$. Since $A \supseteq B$ implies $A \preceq^\triup B$ the set relation $\preceq^\triup$ maintains flexibility: it would be better to have more options than less.

Note that a minimal function value $f(\bar x)$ with respect to $\preceq^\triup$ can exist even if there is no minimal $\bar z \in f(\bar x)$ with respect to $\preceq$.

\begin{example}
Set $f(x) = (x, x+1) \subseteq \R$ with $X = \{x \in \R \mid x \geq 0\}$. Then $x= 0$ satisfies
\[
\forall x \in X \colon f(0) \preceq^\triup f(x)
\]
but $f(0) = (0, 1)$ does not include a minimal element with respect to $\leq$. This example also shows that it does not necessarily make sense just to look for a minimal element $\bar z$ of $\bigcup_{x \in X} f(x)$ and for $\bar x$ such that $\bar z \in f(\bar x)$.
\end{example}

\subsection{Contradictive preferences}

Again, let $\preceq$ be the preference of DM1 and hence $\preceq_2$ the  inverse of $\preceq$, i.e., $\succeq$.

Once $x \in X$ is chosen by DM1, DM2 would pick a maximal element in $f(x)$ w.r.t. $\preceq$ because of the opposite preference, i.e., 
$\bar z \in f(x)$ satisfying
\[
z \in f(x), \; \bar z \preceq z \; \Rightarrow \; z \preceq \bar z
\]
provided there is such an element in $f(x)$.

Thus, DM1 would prefer $x \in X$ over $x' \in X$ if
\[
\forall z \in f(x), \; \exists z' \in f(x') \colon z \preceq z',
\]
i.e., $f(x) \preceq^\trido f(x')$. DM1 would improve or at least not be worse off by selecting $x$ instead of $x'$.

Therefore, DM1 should look for an alternative $\bar x \in X$ such that there is no $x \in X$ satisfying $f(x) \preceq^\trido f(\bar x)$, but not $f(\bar x) \preceq^\trido f(x)$. Thus, DM1 should look for $\bar x \in X$ such that $f(\bar x)$ is a minimal element with respect to the set relation $\preceq^\trido$ . In this case, $A \subseteq B$ implies $A \preceq^\trido B$ which means that DM1 prefers less options for DM2 to inflict damage to her/him. This can also be understood as a preference for inflexibility of the opponent.

Again, a minimal function value $f(\bar x)$ with respect to $\preceq^\trido$ can exist even if there is no maximal $\bar z \in f(\bar x)$  with respect to $\preceq$.

\subsection{Exemplary situations}

The first situation deals with risk evaluation for multivariate positions. The set $X \subseteq L^1_d(\Omega, \F, P)$ represents a collection of future position consisting of $d$ assets where $L^1_d$ denotes the vector space (of equivalence classes) of integrable random variables with values in $\R^d$ over the probability space $(\Omega, \F, P)$: an investor can decide now to invest in these assets and gets a portfolio $x \in X$ at some future point in time. A function $R \colon X \to \P(\R^m)$ is used to evaluate the risk of the positions in terms of $m$ eligible portfolios, e.g., currencies. Usually, one has $1 \leq m \leqq d$. In financial applications of this type, the preferences of the decision makers are often generated by solvency cones modelling transaction costs/bid-ask price spreads between the assets in the portfolio. Background and details can be found, e.g., in \cite{HamelHeyde21Math}. 

The collaborative decision would be if an investor selects $\bar x \in X$ with minimal risk measured by $R(x)$ and then (only later, but maybe by the controlling department of the same investment company) an element $\bar z \in \MIN R(\bar x)$ to be deposited or given as risk compensation. 

A contradictive decision would happen if the investor chooses $\bar x \in X$ and then an external regulator or a bank asks for a maximal (to be on the safe side) deposit out of $R(\bar x)$.

The second situation is a robust optimization model. A function $f \colon X \to \P(Z)$ is generated via $f(x) = \cb{F(x, u) \mid u \in \mathcal U}$ with $F \colon X \times \mathcal U \to Z$ where the set $\mathcal U$ represents a collection of unknown parameters or uncertain data. The decision maker wants to minimize $F$ but does not know which $u \in \mathcal U$ will realize.

The good case, i.e., with aligned preferences, would be if the choice $u \in \mathcal U$ is made in accordance with the preference of the decision maker, i.e., one can minimize $f$ to get $\bar x$ and then a minimal element in $f(\bar x)$ is chosen via selecting $\bar u \in \mathcal U$. In the literature \cite{EhrgottEtAl14EJOR, IdeEtAl14FPTA}, this called the optimistic approach.  

The bad case, i.e., contradicting preferences, would occur if the choice $u \in \mathcal U$ runs against the preference of the decision maker, i.e., it is chosen via maximizing $f(\bar x, u)$ over $u \in \mathcal U$ once $\bar x \in X$ was selected. This is called the pessimistic approach.

Thirdly, a situation arising in multi-stage stochastic programming can be found in \cite{HamelLoehne24ArX}.

\section{Set optimization for decision making}
\label{SecSetOptDM}

So far, the image set for $f$ is just the power set $\P(Z)$. One can order $\P(Z)$ via $\supseteq$ or $\subseteq$ and gets a complete lattice in each case, i.e., every subset has an infimum and a supremum. In particular, $(\P(Z), \supseteq)$ is a complete lattice with 
\begin{equation}
\label{EqPowerInfSupUp}
\inf_{A \in\A} A = \bigcup_{A \in \A} A \quad \text{and} \quad \sup_{A \in \A} A = \bigcap_{A \in \A} A
\end{equation}
for $\A \subseteq \P(Z)$ with the conventions $\inf_{A \in \vn} A = \vn$, $\sup_{A \in \vn} A = Z$. Thus, the least (bottom) element in $(\P(Z), \supseteq)$ is $Z$, the greatest (top) is $\vn$.

Likewise, $(\P(Z), \subseteq)$ is a complete lattice with 
\begin{equation}
\label{EqPowerInfSupDown}
\inf_{B \in \B} B = \bigcap_{B \in \B} B \quad \text{and} \quad \sup_{B \in \B} B = \bigcup_{B \in \B} B
\end{equation}
for $\B \subseteq \P(Z)$ with the conventions $\inf_{B \in \vn} B = Z$, $\sup_{B \in \vn} B = \vn$ with least element $\vn$ and top element $Z$.

Next, two complete lattices are constructed via standard mathematical procedures which correspond to the two set relations defined in the previous section. The remarkable fact is that the infimum and supremum in these lattices are meaningful concepts--even though minimal or maximal elements do not need to exist, neither in $\P(Z)$ with respect to $\preceq^\triup$ or $\preceq^\trido$, nor in $Z$ with respect to $\preceq_1$ or $\preceq_2$.

As a result, the first stage decision making problem becomes an optimization problem with a complete lattice-valued objective function for which sound solution concepts can be introduced.

\subsection{The collaborative case}
\label{SubSecCollaborative}

The construction of the complete lattice is as follows. For $A \subseteq Z$, define
\[
A^\triup = \cb{z \in Z \mid \exists a \in A \colon a \preceq z}
\]
with $\vn^\triup = \vn$.

\begin{lemma}
\label{LemUpperLattice}
The function $A \mapsto A^\triup$ is a closure operator, i.e., it satisfies

(C1u) $A^\triup \supseteq A$,

(C2u) $(A^\triup)^\triup = A^\triup$ (idempotency),

(C3u) $(A \cup B)^\triup = A^\triup \cup B^\triup$ (preservation of binary unions)
\\
for $A \in \P(Z)$. The set
\[
\P^\triup(Z, \preceq) = \cb{A \subseteq Z \mid A^\triup = A}
\]
is a complete lattice with respect to $\supseteq$, and one has for all $A, B \in \P^\triup(Z, \preceq)$
\[
A \supseteq B \quad \Leftrightarrow \quad A \preceq^\triup B.
\]
\end{lemma}

It is well-known that property (C3u) readily implies monotonicity, i.e., $^\triup$ is order preserving: $B \subseteq A$ implies $B^\triup \subseteq A^\triup$: indeed, since $A \cup B = A$ in this case, (C3u) implies $A^\triup = A^\triup \cup B^\triup \supseteq B^\triup$. Thus, the $^\triup$-operation is an extensive Kuratowski closure operator (\cite{Kuratowski22FM}): it satisfies (C1u)-(C3u) together with $\vn^\triup = \vn$.

{\sc Proof of Lemma \ref{LemUpperLattice}.} (C1u) and (C2u) are straightforward. (C3u): Assume $z \in (A \cup B)^\triup$. Then, there is $z' \in A \cup B$ with $z \preceq z$, hence $z \in A^\triup \cup B^\triup$. Conversely, if $z \in A^\triup \cup B^\triup$, then there is $a \in A$ with $a \preceq z$ or $b \in B$ with $b \preceq z$. Thus $z \in (A \cup B)^\triup$. 

For the equivalence of the orders relations, one has to verify that $A \preceq^\triup B$ implies $A \supseteq B$ for $A, B \in \P^\triup(Z, \preceq)$. Indeed, $A \preceq^\triup B$ means that for each $b \in B$ there is $a \in A$ with $a \preceq b$. Since $A \in \P^\triup(Z, \preceq)$, this implies $b \in A$, thus $B \subseteq A$.

Finally, one can apply the Knaster-Tarski fixed point theorem to the complete lattice $(\P(Z), \supseteq)$ and the closure operator $^\triup$ which is monotone (order preserving). Thus, the set of fixed points of the function $A \mapsto A^\triup$, namely $\P^\triup(Z, \preceq)$, is a complete lattice with respect to $\supseteq$. \pend

One may note that major parts of the above result are well-known in lattice theory; see \cite{DaveyPriestley02Book2nd} from which one can combine 7.3 Theorem, 7.5 Example (5) and 1.27/1.28 to get the complete lattice structure of $\P^\triup(Z, \preceq)$. However, no link to set relations is given therein and neither is (C3u).

\begin{remark}
\label{RemUCSs} 
The set $\P^\triup(Z, \preceq)$ can be obtained in a different way. The equivalence relation $S(R^\triup)$ is given by $(A, B) \in S(R^\triup)$ if, and only if,   $A \preceq^\triup B$, $B \preceq^\triup A$. The set $A^\triup$ is the largest set in the equivalence class with respect to $S(R^\triup)$ generated by $A$:
\[
A^\triup = \bigcup \cb{B \subseteq Z \mid (A, B) \in S(R^\triup)}.
\]
It can thus be understood as the minimal representative of the equivalence class generated by $A$ since the union is the infimum in $(\P(Z), \supseteq)$.
\end{remark}

\begin{proposition}
\label{PropTriupUnion}
The closure operator $^\triup$ commutes with the infimum, i.e.,
\[
\of{\inf_{A \in \A} A}^\triup = \inf_{A \in \A} A^\triup
\]
while it does not commute with the supremum in general.
\end{proposition}

{\sc Proof.} Since the infimum is the union as given in \eqref{EqPowerInfSupUp} one has to verify $\of{\bigcup_{A \in \A} A}^\triup = \bigcup_{A \in \A} A^\triup$ which can easily be done. A counterexample for the supremum case is the following: take $Z = \R^2$, $\preceq$ as the component-wise order and consider $\mathcal A = \{(1,0)^\top, (0,1)^\top\}$ with $A_1 = \{(1,0)^\top\}$, $A_2 = \{(0, 1)^\top\}$. Then $\sup\{A_1, A_2\} = \vn$ while $(A_1)^\triup = \{(1,0)^\top\} + \R^2_+$, $(A_2)^\triup = \{(0, 1)^\top\} + \R^2_+$ and $\sup\{(A_1)^\triup, (A_2)^\triup\} =   \{(1, 1)^\top\} + \R^2_+ \neq \vn$. \pend

\begin{lemma}
\label{LemUpperInfMin}
(1) For $A \in \P(Z)$ one has
\[
z \in \MIN A^\triup \; \Leftrightarrow \; \exists a \in \MIN A \colon a \preceq z \preceq a.
\]
If $\preceq$ is additionally antisymmetric, then $\MIN A = \MIN A^\triup$.

(2) For $\A \subseteq  \P^\triup(Z, \preceq)$ one has
\[
\inf_{A \in \A} A = \bigcup_{A \in \A} A \quad \text{and} \quad \sup_{A \in \A} A = \bigcap_{A \in \A} A.
\]
\end{lemma}

{\sc Proof.} (1) If $z \in \MIN A^\triup$, then $z \in A^\triup$ and there is $a \in A$ with $a \preceq z$. From $A \subseteq A^\triup$ and $z \in \MIN A^\triup$ one gets $z \preceq a$. Vice versa, if $a \in \MIN A$, $z' \in A^\triup$ with $a \preceq z \preceq a$ and $z' \preceq z$, then $z' \preceq a$ and thus $z' \preceq z$ both by transitivity. Hence $z \in  \MIN A^\triup$. The last assertion is immediate from antisymmetry.

(2) The formulas are consequences of \eqref{EqPowerInfSupUp} and Lemma \ref{LemUpperLattice}. \pend

Turning back to the decision making problem, one may note that the decisions are not altered if $f$ is replaced by $f^\triup$ with $f^\triup(x) = (f(x))^\triup$: this is a consequence of Lemma \ref{LemUpperLattice} and Lemma \ref{LemUpperInfMin}. The infimum
\[
\inf_{x \in X} f^\triup(x) = \bigcup_{x \in X} f^\triup(x)
\]
collects all outcomes of the second stage decision along with alternatives which are worse than potential outcomes.

The first stage decision now becomes a set optimization problem with an objective function mapping into a complete lattice of sets:
\[
\tag{SOColl} \text{minimize} \quad f^\triup(x) \quad \text{subject to} \quad x \in X.
\]
The meaning of "minimize" is to look for the infimum and minimal function values.

\begin{definition}
\label{DefSolConceptColl}
A set $S \subseteq X$ is called a solution of (SOColl) if

(S1) $\inf_{y \in S} f^\triup(y) = \inf_{x \in X} f^\triup(x)$ and

(S2) $f^\triup(y)$ is minimal in $\cb{f^\triup(x) \mid x \in X}$ with respect to $\supseteq$ for all $y \in S$.
\end{definition}

This solution concept is due to \cite{HeydeLoehne11Opt} with precursor \cite{Hamel04Proc}. An element $\bar x \in X$ such that $f^\triup(\bar x)$ is minimal in $\cb{f^\triup(x) \mid x \in X}$ is called a minimizer of $f^\triup$. Thus, Definition \ref{DefSolConceptColl} asks for a set $S \subseteq X$ which produces the infimum and only includes minimizers.

\begin{example}
\label{ExSolutionIssues}
This example is \cite[Example 1.2]{Loehne24ArX}. Let three subsets of $Z = \R^2$ be defined by $Z_1 = \{(1,0)^\top\} + \R^2_+$, $Z_2 = \{(0,1)^\top\} + \R^2_+$, $Z_3 = \co\{(1.05,0.05)^\top, (0.05,1.05)^\top\} + \R^2_+$. The function $f \colon \R^3 \to \P(\R^2)$ defined by
\[
f(x) = \left\{
	\begin{array}{ccc}
	x_1 Z_1 + x_2 Z_2 + x_3Z_3 & \colon & x_1, x_2, x_3 \geq 0, x_1 + x_2 + x_3 = 1 \\
	\vn & \colon & \text{otherwise}
	\end{array}
	\right.
\]
actually maps into $\P^\triup(\R^2, \leq_{\R^2_+})$ which means $f = f^\triup$. The set $S = \{x \in \R^3 \mid x_1, x_2 \geq 0, x_3 = 0, x_1 + x_2 = 1\}$ is a solution of (SOColl). However, its elements do not provide much flexibility for the decision maker since the function values $f(x) = \{(x_1, x_2)^\top\} + \R^2_+$ are all minimal with respect to $\supseteq$ but have the only minimal element $(x_1, x_2)$. The function value $f(\bar x) = Z_3$ for $\bar x = (0,0,1)^\top$ is also minimal and provides a much bigger set of minimal elements. However, $\bar x \not\in S$. On the other hand, $\bar x$ can be included in any solution set since $f(\bar x)$ is minimal.
\end{example}

A decision procedure for the collaborative problem is as follows.

\medskip\noindent
{\bf Decision Procedure DPColl1.}
\begin{description}
\item[{\bf 1}] Determine a solution $S \subseteq X$ of (SOColl) and select $\bar x \in S$.
\item[{\bf 2}] Select $\bar z \in \MIN f(\bar x)$. 
\end{description}

A naive application of DPColl1 might run into a number of issues. 

The first is that a solution to (SOColl1) might not exist or is not unique (see the example above). While in the first case, one can resort to approximate solutions (see \cite{CrespiEtAl21MOR}), there is no criterion yet how to choose one among several solutions for (SCOll): one might prefer solutions whose elements provide greater flexibility. One might be tempted to say that the same issue arises for scalar optimization problems with several solutions; however, the function values can be very different in the set-valued case while in the scalar case they all coincide.

The second issue is that $\MIN f(\bar x)$ might be empty for the selected $\bar x$ or not to the liking of the DM for reasons outside the mathematical model.

For both issues, conditions can be found excluding them such as a Weierstrass type existence theorem for (SOColl) along with lower domination assumptions to the values of $f$, see \cite{HeydeLoehne11Opt}, \cite{CrespiEtAl21MOR}. This will not be discussed here. Instead, a second decision making procedure is outlined which also addresses the second issue mentioned above. The idea is due to \cite{Loehne24ArX} and was applied to multi-stage stochastic optimization problems in \cite{HamelLoehne24ArX}.

\medskip
\noindent
{\bf Decision Procedure DPColl2.}
\begin{description}
\item[{\bf 1}] Select a desirable outcome set $Y \subseteq Z$, for example $Y = \{y_1, \ldots,y_N\}$ with a positive integer $N$.
\item[{\bf 2}] Determine $\bar x \in X$ such that $Y \subseteq f^\triup(\bar x)$ and $f^\triup(\bar x)$ is minimal in $\cb{f^\triup(x) \mid x \in X}$ with respect to $\supseteq$.
\item[{\bf 3}] Select $\bar z \in \MIN f^\triup(\bar x)$.
\end{description}

Again, there might not even be an $x \in X$ with $Y \subseteq f^\triup(x)$ in which case DM has to modify $Y$ in order to align her/his wishes with the reality (of the mathematical model). This leads to a trial-and-error process which turns out to be finite in the polyhedral convex case as shown in \cite{Loehne24ArX}. Note that $\bar z$ in Step 3 does not need to be an element of $Y$.

DPColl1 can be considered as a forward procedure while DPColl2 starts from desirable outcomes and thus works backwards in the sense that the first stage decision is directed by second stage desirable outcomes. The following result shows that both procedures lead to the same type of problem.

\begin{theorem}
\label{ThmDesirableOutcomesColl}
Let $Y \subseteq Z$ and define $\P^\triup_Y(Z, \preceq) = \cb{A \in \P^\triup(Z, \preceq) \mid Y \subseteq A} \cup\{\vn\}$. Then $(\P^\triup_Y(Z, \preceq), \supseteq)$ is a sublattice of $(\P^\triup(Z, \preceq), \supseteq)$ with greatest element $\vn$ and least element $Z$.
\end{theorem}

{\sc Proof.} Clearly, $\supseteq$ is a partial order on $\P^\triup_Y(Z, \preceq)$. It remains to show that 
\[
\bigcup_{A \in \mathcal A} A \quad\text{and}\quad \bigcap_{A \in \mathcal A} A
\]
are in $\P^\triup_Y(Z, \preceq)$ whenever $\mathcal A \subseteq \P^\triup_Y(Z, \preceq)$ which is obvious (recall the conventions for $\mathcal A = \vn$ from \eqref{EqPowerInfSupUp}). Moreover, $Z \supseteq A \supseteq Y^\triup \supseteq \vn$ is clearly true for all $A \in \P^\triup_Y(Z, \preceq)$. \pend

As a consequence, one can re-define $f^\triup$ by setting $f^\triup_Y(x) = f^\triup(x)$ if $Y \subseteq f^\triup(x)$, $f^\triup_Y(x) = \vn$ if $Y \not\subseteq f^\triup(x)$ and solve (SOColl) with $f^\triup$ replaced by $f^\triup_Y$. If $\bar x$ belongs to a solution set for the redefined problem, then $Y \subseteq f^\triup(x)$ is automatically satisfied. Thus, DPColl2 can be re-phrased as: apply DPColl1 to $f^\triup_Y$ instead of $f^\triup$. This leaves the question for the relationship between solutions of the set minimization problems for $f^\triup$ and $f^\triup_Y$. 

The domain of the function $f^\triup_Y$ is the set
\[
\dom f^\triup_Y = \cb{x \in X \mid f^\triup_Y(x) \neq \vn}
\]
with $\vn$ being the greatest element in $(P^\triup_Y(Z, \preceq), \supseteq)$. One may observe that a reasonable requirement to $Y$ is the condition
\begin{equation}
\label{EqConditionToYup}
\inf_{x \in X} f^\triup(x) = \bigcup_{x \in X} f^\triup(x) \supseteq Y^\triup.
\end{equation}
Of course, it does not guarantee $\vn \neq \dom f^\triup_Y$, but if it is not satisfied, then $\dom f^\triup_Y = \vn$ for sure. Because $\supseteq$ is the order in $P^\triup_Y(Z, \preceq)$, \eqref{EqConditionToYup} translates into ''the infimum of $f^\triup$ is less than or equal to $Y^\triup$." Thus, $Y^\triup$ shouldn't include "too many" desirable outcomes.

In the following, it is shown that results from \cite{Loehne24ArX} carry over to the much more general case considered in this paper.

Define the function $v^\triup_f \colon \P(Z) \to \P^\triup(Z, \preceq)$ via
\[
v^\triup_f(Y) = \inf_{x \in X} f^\triup_Y(x) = \bigcup_{x \in X} f^\triup_Y(x) = \bigcup_{x \in X}\cb{f^\triup(x) \mid x \in X, \, Y \subseteq f^\triup(x)}.
\]
Note that $f^\triup_Y$ does not appear in \cite[p. 7]{Loehne24ArX}. Moreover, one can assume $Y \in  \P^\triup(Z, \preceq)$ without loss of generality since $f^\triup_Y(x) \in \P^\triup(Z, \preceq)$.

\begin{proposition}
\label{PropValueFunctionProp}
(1) If there is $x \in X$ with $Y \subseteq f^\triup(x)$, then $Y \subseteq v^\triup_f(Y)$ for $Y \in \P^\triup(Z, \preceq)$.

(2) One has $v^\triup_f(\vn) = \inf_{x \in X} f^\triup(x)$.

(3) $Y, Y' \in \P^\triup(Z, \preceq)$ and $Y \subseteq Y'$ imply $v^\triup_f(Y') \subseteq v^\triup_f(Y)$.
\end{proposition}

{\sc Proof.} (1) This follows from the definitions of $v^\triup_f(Y)$ and $f^\triup_Y(x)$ since $Y \subseteq f^\triup_Y(x)$ is true for all $x \in X$.

(2) It holds
\[
v^\triup_f(\vn) = \cb{z \in Z \mid \exists x \in X \colon z \in f^\triup_\vn(x)}.
\]
The result follows since $f^\triup_\vn(x) = f^\triup(x)$ for all $x \in X$ by definition.

(3) If $v^\triup_f(Y') = \vn$, then the inclusion $v^\triup_f(Y') \subseteq v^\triup_f(Y)$ trivial. If $z \in v^\triup_f(Y')$, then there is $x \in X$ such that $z \in f^\triup(x)$ and $Y' \subseteq f^\triup(x)$. By assumption, $Y \subseteq f^\triup(x)$ and hence $z \in f^\triup_Y(x)$ which turn gives $z \in v^\triup_f(Y)$.
\pend

It turns out that minimizers of $f^\triup$ can be characterized via equality in (1) of the previous proposition.

\begin{definition}
\label{DefFixedPointValue}
An element $Y \in \P^\triup(Z, \preceq)$ is called a fixed point of a function $v \colon \P^\triup(Z, \preceq) \to \P^\triup(Z, \preceq)$ if
$Y = v(Y)$.
\end{definition}

The following result is a general version of \cite[Proposition 3.1]{Loehne24ArX} which deals with the case of $X$ and $Z$ being finite dimensional vector spaces.

\begin{theorem}
\label{ThmFixedPointValue}
The following statements are equivalent for $Y \in \P^\triup(Z, \preceq)$:

(1) $Y$ is a fixed point of $v^\triup_f$.

(2) $f^\triup(x) = Y$ for all $x \in \dom f^\triup_Y$.

(3) There is a minimizer $\bar x \in X$ of $f^\triup$ with $Y = f^\triup(\bar x)$.
\end{theorem}

{\sc Proof.} (1) $\Rightarrow$ (2). If $Y = \vn$, then $v^\triup_f(Y) = \cup_{x \in X}f^\triup(x)$, thus $v^\triup_f(\vn) = \vn$ if, and only if, $f^\triup(x) = \vn = Y$ for all $x \in X$, so (2) is satisfied. If $Y \neq \vn$ there must be $x \in X$ with $Y \subseteq f^\triup(x) \neq \vn$ since $Y$ is a fixed point of $v^\triup_f$. If $z \in f^\triup(x)\bs Y$ for such $x$, then $z \in v^\triup_f(Y)\bs Y$, a contradiction. Thus, (2) is satisfied also in this case.

(2) $\Rightarrow$ (3). If $Y = \vn$, then $f^\triup(x) = \vn$ for all $x \in X$, thus every $x \in X$ is a minimizer with $Y = f^\triup(x)$. Take $\bar x \in \dom f^\triup_Y$, i.e., $Y \subseteq f^\triup(\bar x) = Y$. If $f^\triup(x) \supseteq f^\triup(\bar x) = Y$, then $x \in \dom f^\triup_Y$, hence $f^\triup(x) = Y = f^\triup(x)$, thus $\bar x$ is a minimizer of $f^\triup$. 

(3) $\Rightarrow$ (1). Let $\bar x \in X$ be a minimizer of $f^\triup$ with $Y = f^\triup(\bar x)$. If $\dom f = \vn$, then $Y = \vn$ and $v^\triup_f(\vn) = \vn$. If $\dom f \neq \vn$, then $Y = f^\triup(\bar x) \neq \vn$ and there is no $x \in X$ with $f^\triup(x) \supsetneq f^\triup(\bar x)$, hence $Y \subseteq v^\triup_f(Y)$ by definition of $v^\triup_f$. If $Y \neq v^\triup_f(Y)$ would be true, then there would be $z \in  v^\triup_f(Y)\bs Y$ and hence $x \in X$ with $z \in f^\triup_Y(x)$ by definition of $v^\triup_f$. This would imply $f^\triup_Y(x) \supsetneq f^\triup_Y(\bar x)$ and consequently $f^\triup(x) \supsetneq f^\triup(\bar x)$ contradicting the minimality of $\bar x$. \pend

Proposition \ref{PropValueFunctionProp} and Theorem \ref{ThmFixedPointValue} provide an optimality condition for Problem (SOColl): a solution is a set $S$ satisfying $\inf_{x \in S} f^\triup(x) = \inf_{x \in X} f^\triup(x)$ and $f^\triup(x)$ is a fixed point for $v^\triup_f$ for all $x \in S$.

\begin{remark}
\label{RemMinMin}
In the collaborative case, the decision making problem is of the "min-min" type, i.e., a relatively easy one from a mathematical point of view reflecting the fact that collaboration under aligned preferences is usually easier than solving a conflict situation.
\end{remark}

\subsection{The contradictive case}
\label{SubSecConflict}

The construction of the complete lattice is as follows. For $A \subseteq Z$, define
\[
A^\trido = \cb{z \in Z \mid \exists a \in A \colon z \preceq a}
\]
with $\vn^\trido = \vn$.

\begin{lemma}
\label{LemLowerLattice}
The function $A \mapsto A^\trido$ is a closure operator, i.e., it satisfies

(C1l) $A^\trido \supseteq A$,

(C2l) $(A^\trido)^\trido = A^\trido$ (idempotency),

(C3l) $(A \cup B)^\trido = A^\trido \cup B^\trido$ (preservation of binary unions)
\\
for $A \in \P(Z)$. The set
\[
\P^\trido(Z, \preceq) = \cb{A \subset Z \mid A^\trido = A}
\]
is a complete lattice with respect to $\subseteq$, and one has for all $A, B \in \P^\trido(Z, \preceq)$
\[
A \subseteq B \quad \Leftrightarrow \quad A \preceq^\trido B.
\]
\end{lemma}

{\sc Proof.} If one defines a new preorder $\preceq'$ on $Z$ by $y \preceq' z$ if $z \preceq y$, then one may see that 
$A^\trido = A^{\triup'}$ which makes the current lemma a version of Lemma \ref{LemLowerLattice}. In particular, the latter gives $A \subseteq B$ if, and only if, $B \preceq ^{\triup'} A$ which in turn is $A \preceq^\trido B$ for $A, B \in \P^\trido(Z, \preceq)$. \pend

Clearly, (C3l) also implies monotonicity: $B \subseteq A$ implies $B^\trido \subseteq A^\trido$, and the $^\trido$-operation is an extensive Kuratowski closure operator. Moreover, Proposition \ref{PropTriupUnion} remains valid in this case, i.e., $^\trido$ commutes with the supremum which is a union in this case according to \eqref{EqPowerInfSupDown}, but it does not commute with the infimum.

\begin{lemma}
\label{LemLowerSupMax}
(1) For $B \in \P(Z)$ one has
\[
z \in \MAX B^\trido \; \Leftrightarrow \; \exists b \in \MAX B \colon b \preceq z \preceq b.
\]
If $\preceq$ is additionally antisymmetric, then $\MAX B = \MAX B^\trido$.

(2) For $\B \subseteq  \P^\trido(Z, \preceq)$ one has
\[
\inf_{B \in \B} B = \bigcap_{B \in \B} B \quad \text{and} \quad \sup_{B \in \B} B = \bigcup_{B \in \A} B.
\]
\end{lemma}

{\sc Proof.} Again, by replacing $\preceq$ by $\preceq'$ as in the proof of Lemma \ref{LemLowerLattice} one can reduce the result to Lemma \ref{LemUpperInfMin}. \pend

Thus, the decision making problem is not altered if $f$ is replaced by $f^\trido$ with $f^\trido(x) = (f(x))^\trido$: one only adds elements to $f(x)$ which would not be chosen by DM2 anyway. The supremum
\[
\sup_{x \in X} f^\trido(x) = \bigcup_{x \in X} f^\trido(x)
\]
collects all outcomes of the second stage decision along with alternatives which are better for DM1 than potential outcomes. It would be good for DM1 to know this set if she/he is interested in a worst case analysis. On the other hand, it might seem more important for DM1 to know the unavoidable outcomes since DM2 is going to "punish" her/him. These outcomes are collected in the set
\[
\inf_{x \in X} f^\trido(x) = \bigcap_{x \in X} f^\trido(x).
\]
Indeed, whatever DM1 decides, elements of this set are available to DM2. However, it can be empty  and even if not, it can be "far away" from any alternative which will be available for DM2.

\begin{example}
\label{ExIntersectionIssue}
(1) Let $X = \{x \in \R^2_+ \mid x_1+x_2 = s\}$ for fixed $s>0$ and $f \colon X \to \P(\R^2)$ be defined by $f(x) = \{x\}$. The preference is the component-wise order generated by the convex cone $\R^2_+$, denoted $\leq_{\R^2_+}$. Then $f^\trido(x) = \{x\} - \R^2_+$ and
\[
        \inf_{x \in X} f^\trido(x) = \bigcap_{x \in X} f^\trido(x) = \{(0,0)^\top\} -\R^2_+.
\]
The distance between the infimum and the set $\{f(x) \mid x \in X\}$ can be made arbitrary big by increasing $s$.

(2) Let $X = \{x \in \R^2 \mid 2x_1+x_2 \leq 4, x_1 + 2x_2 \leq 4\}$ and $f(x) = \{x\}$ with the same preference as in (1). Then $\bigcap_{x \in X} f^\trido(x) = \vn$.
\end{example}

On the other hand, DM1 should select $\bar x \in X$ such that $f^\trido(\bar x)$ is minimal, i.e.,
\[
x \in \dom f, f^\trido(\bar x) \subseteq f^\trido(x) \quad \Rightarrow \quad f^\trido(x) = f^\trido(\bar x).   
\]

The first stage decision again becomes a set optimization problem with an objective function mapping into a complete lattice of sets:
\[
\tag{SOCon} \text{minimize} \quad f^\trido(x) \quad \text{subject to} \quad x \in X.
\]
The meaning of "minimize" is to look for the infimum and minimal function values, but it now constitutes a different mathematical problem since the infimum is an intersection instead of a union.

\begin{definition}
\label{DefSolConceptConf}
A set $T \subseteq X$ is called a solution of (SOCon) if

(T1) $\inf_{y \in T} f^\trido(y) = \inf_{x \in X} f^\trido(x)$ and

(T2) $f^\trido(y)$ is minimal in $\cb{f^\trido(x) \mid x \in X}$ with respect to $\subseteq$ for all $y \in T$.
\end{definition}

As in the collaborative case, an element $\bar x \in X$ such that $f^\trido(\bar x)$ is minimal in $\cb{f^\trido(x) \mid x \in X}$ is called a minimizer of $f^\trido$. The decision procedure now looks as follows.

\medskip
\noindent
{\bf Decision Procedure DPCon1.}
\begin{description}
\item[{\bf 1}] DM1 determines a solution $T \subseteq X$ of (SOCon) and selects $\bar x \in T$.
\item[{\bf 2}] DM2 selects $\bar z \in \MAX f(\bar x)$. 
\end{description}

Definition \ref{DefSolConceptConf} and DPCon1 can both be considered as a naive translation of Definition \ref{DefSolConceptColl} and DPColl1 to the contradictive case. The issues discussed above and Example \ref{ExIntersectionIssue} show that it might not be the best approach. Moreover, as in the collaborative case, DPCon1 is not guided by potential outcomes. An analog to DPColl2 is as follows.

\medskip
\noindent
{\bf Decision Procedure DPCon2.}
\begin{description}
\item[{\bf 1}] Select a set $Y \subseteq Z$ of outcomes which are acceptable for DM1, for example $Y = \{y_1, \ldots,y_N\}$ with a positive integer $N$.
\item[{\bf 2}] Determine $\bar x \in X$ such that $f^\trido(\bar x) \subseteq Y^\trido$ and $f^\trido(\bar x)$ is minimal in $\cb{f^\trido(x) \mid x \in X}$ with respect to $\subseteq$.
\item[{\bf 3}] DM2 selects $\bar z \in \MAX f^\trido(\bar x)$.
\end{description}

Since the set $Y$ includes acceptable outcomes, DM1 looks for $x$ such that $f^\trido(x) \subseteq Y^\trido$ which means that $f^\trido(x)$ would be ''more than acceptable." The next theorem corresponds to Theorem \ref{ThmDesirableOutcomesColl}.

\begin{theorem}
\label{ThmDesirableOutcomesConf}
Let $Y \subseteq Z$ and define $\P^\trido_Y(Z, \preceq) = \cb{B \in \P^\trido(Z, \preceq) \mid B \subseteq Y} \cup\{Z\}$. Then $(\P^\trido_Y(Z, \preceq), \subseteq)$ is a sublattice of $(\P^\trido(Z, \preceq), \subseteq)$ with greatest element $Z$ and least element $\vn$.
\end{theorem}

{\sc Proof.} Clearly, $\subseteq$ is a partial order on $\P^\trido_Y(Z, \preceq)$. One has that
\[
\bigcup_{B \in \mathcal B} B \quad\text{and}\quad \bigcap_{B \in \mathcal B} B
\]
are in $\P^\trido_Y(Z, \preceq)$ whenever $\mathcal B \subseteq \P^\trido_Y(Z, \preceq)$. Clearly, $\vn$ is the least element in  $\P^\trido_Y(Z, \preceq)$ and $Z$ is the greatest since the order relation is $\subseteq$. \pend

One can now re-define $f^\trido$ by setting $f^\trido_Y(x) = f^\trido(x)$ if $f^\trido(x) \subseteq Y$, $f^\trido_Y(x) = Z$ if $ f^\trido (x) \not\subseteq Y$ and solve (SOCon) with $f^\trido$ replaced by $f^\trido_Y$. The domain of the function $f^\trido_Y$ is the set
\[
\dom f^\trido_Y = \cb{x \in X \mid f^\trido_Y(x) \neq Z}
\]
recalling that $Z$ is the greatest element in $(P^\trido_Y(Z, \preceq), \subseteq)$. 
Parallel to the collaborative case, a reasonable requirement to $Y$ is the condition
\begin{equation}
\label{EqConditionToYdown}
\inf_{x \in X} f^\trido(x) = \bigcap_{x \in X} f^\trido(x) \subseteq Y^\trido
\end{equation}
which does not guarantee $\vn \neq \dom f^\trido_Y$, but if it is not satisfied, then $\dom f^\trido_Y = \vn$ for sure. Because $\subseteq$ is the order in $P^\trido_Y(Z, \preceq)$, \eqref{EqConditionToYdown} translates into ''the infimum of $f^\trido$ is less than or equal to $Y^\trido$." Thus, $Y^\trido$ should include "enough" acceptable outcomes.

If $\bar x$ belongs to a solution set for the re-defined problem, then $f^\trido(x) \subseteq Y$ is automatically satisfied. Thus, DPCon2 can be re-phrased as: apply DPCon1 to $f^\trido_Y$ instead of $f^\trido$. The leaves the question for the relationship between solutions of the set minimization problems for $f^\triup$ and $f^\triup_Y$. 

Define the function $v^\trido_f \colon \P(Z) \to \P^\trido(Z, \preceq)$ via
\[
v^\trido_f(Y) = \bigcap_{x \in X} f^\trido_Y(x) = \inf_{x \in X} f^\trido_Y(x) = \bigcap\cb{f^\trido (x) \mid x \in X, \, f^\trido(x) \subseteq Y}.
\]
Moreover, one can assume $Y \in  \P^\trido(Z, \preceq)$ without loss of generality since $f^\trido_Y(x) \in \P^\trido(Z, \preceq)$.

\begin{proposition}
\label{PropValueFunctionProp}
(1) If there is $x \in X$ with $f^\trido(x) \subseteq Y$, then $v^\trido_f(Y) \subseteq Y$ for $Y \in \P^\trido(Z, \preceq)$.

(2) One has $v^\trido_f(Z) = \inf_{x \in X} f^\trido(x)$.

(3) $Y, Y' \in \P^\trido(Z, \preceq)$ and $Y \subseteq Y'$ imply $v^\trido_f(Y') \subseteq v^\trido_f(Y)$.
\end{proposition}

{\sc Proof.} (1) This follows from the definition of $v^\trido_f(Y)$ since $f^\trido(x) \subseteq Y$ is true for at least one $x \in X$.

(2) It holds
\[
v^\trido_f(Z) = \bigcap\cb{f^\trido (x) \mid x \in X, \, f^\trido(x) \subseteq Z} = \inf_{x \in X} f^\trido(x).
\]

(3) The assumption $Y \subseteq Y'$ implies $\{x \in X \mid f^\trido(x) \subseteq Y\} \subseteq \{x \in X \mid f^\trido(x) \subseteq Y'\}$. Since $v^\trido_f$ is an intersection, one gets  $v^\trido_f(Y') \subseteq v^\trido_f(Y)$.
\pend

It turns out that minimizers of $f^\triup$ can be characterized via equality in (1) of the previous proposition.

\begin{definition}
\label{DefFixedPointValueDo}
An element $Y \in \P^\trido(Z, \preceq)$ is called a fixed point of a function $v \colon \P^\trido(Z, \preceq) \to \P^\trido(Z, \preceq)$ if
$Y = v(Y)$.
\end{definition}

The following result parallels Theorem \ref{ThmFixedPointValue} and is new. 

\begin{theorem}
\label{ThmFixedPointValueDo}
The following statements are equivalent for $Y \in \P^\trido(Z, \preceq)$:

(1) $Y$ is a fixed point of $v^\trido_f$.

(2) $f^\trido(x) = Y$ for all $x \in \dom f^\trido_Y$.

(3) There is a minimizer $\bar x \in X$ of $f^\trido$ with $Y = f^\trido(\bar x)$.
\end{theorem}

{\sc Proof.} (1) $\Rightarrow$ (2). If $Y = Z$, then $v^\triup_f(Y) = \cap_{x \in X}f^\trido(x) = Z$, thus $f^\trido(x) = Z$ for all $x \in X$, so (2) is satisfied. If $Y \neq Z$ there must be $x \in X$ with $Z \neq f^\trido(x)  \subseteq Y$ since $Y$ is a fixed point of $v^\trido_f$. Take such an $x$ and assume $z \in Y\bs f^\trido(x)$. Then $z \notin \bigcap\cb{f^\trido (x) \mid x \in X, \, f^\trido(x) \subseteq Y}$ contradicting $Y = v^\trido_f(Y)$. Thus, (2) is satisfied also in this case.

(2) $\Rightarrow$ (3). If $Y = Z$, then $f^\triup(x) = Z$ for all $x \in X$, thus every $x \in X$ is a minimizer with $Y = Z = f^\trido(x)$. If $\bar x \in \dom f^\trido_Y$ satisfies $Y = f^\trido(\bar x)$, then $\bar x$ is a minimizer since there is no $x \in X$ with $f^\trido(x) \subsetneq f^\trido(\bar x) = Y$ by assumption. 

(3) $\Rightarrow$ (1). By definition of $v^\trido_f$,
\[
v^\trido_f(f^\trido ( \bar x)) = \bigcap\cb{f^\trido(x) \mid f^\trido(x) \subseteq f^\trido(\bar x)} = f^\trido(\bar x)
\]
since $\bar x$ is a minimizer of $f^\trido$ with respect to $\subseteq$. Because of the assumption $Y = f^\trido(\bar x)$, (1) is satisfied.

\begin{remark}
\label{RemMinMax}
In the contradictive case, the decision making problem is of the "min-max" type which makes it mathematically more difficult to solve reflecting the fact that a conflict usually is difficult to deal with. 
\end{remark}

\section{Partially aligned preferences}
\label{SecPartiallyAligned}

It might be realistic to assume that the two decision makers do not have completely opposed preferences but could agree on some (non-trivial) comparisons. Thus, $R_1 \cap R_2 \neq D$ which means that there are pairs $y, z \in Z$ with $y \preceq_1 z$, $y \preceq_2 z$ and $y \neq z$. Clearly, $R := R_1 \cap R_2$ is a preorder (a partial order) if $R_1$, $R_2$ are.

In such a case, it might also happen that the two decision makers do not agree and even are willing to give up an advantage in order to "punish" the opponent. Thus, they might consider to make decisions according to $\hat R_1 := \sqb{R_1\bs (R_1 \cap R_2)} \cup D$ and $\hat R_2 := \sqb{R_2\bs (R_1 \cap R_2)} \cup D$, respectively, which can be considered as an extreme behaviour. This can lead to different outcomes.

\begin{example}
\label{ExDM1Dilemma}
Consider the function
\[
f(x)  =
\left\{
\begin{array}{ccc}
\co\{(0, 2x), (1-2x, 1)\} & \colon & 0 \leq x \leq \frac{1}{2} \\
\co\{(2-2x, 0), (1, -1 + 2x)\} & \colon & \frac{1}{2} \leq x \leq 1\\
\vn & \colon & x < 0 \;\, \text{or} \;\, x > 1.
\end{array}
\right.
\]
The image set, i.e., $f[X] = \{z \in \R^2 \mid \exists x \in X \colon z \in f(x)\}$ with $X = [0,1]$ is the square $Q$ with the vertices $(0,0)$, $(1,0)$, $(0,1)$, $(1,1)$ while the function values are segments parallel to the diagonal of $Q$ connecting $(0,0)$ and $(1,1)$.

Assume DM1 has the preference $R_1$ which is the vector order generated by the convex cone $C_1 := \R^2_+$, i.e., the component-wise order. This gives
\[
\inf_{x \in X} f^{\triup_1}(x) = Q + \R^2_+, \quad \inf_{x \in X} f^{\trido_1}(x) = \{(0,0)\} - \R^2_+
\]
where the first infimum is in $\P^{\triup_1}(\R^2, \leq_{\R^2_+})$ with respect to $\supseteq$, i.e., it is a union, while the second infimum is in  $\P^{\trido_1}(\R^2, \leq_{\R^2_+})$ with respect to $\subseteq$, i.e., it is an intersection. The singleton $\{\bar x = \frac{1}{2}\}$ produces the infimum of $f^{\triup_1}$ since $f^{\triup_1}(1/2) = Q + \R^2_+$ while the set $\{0, 1\}$ produces the infimum of $f^{\trido_1}$ since $\{(0,0)\} - \R^2_+ = f^{\trido_1}(0) \cap f^{\trido_1}(1)$. Moreover, $f^{\triup_1}(1/2) \supseteq f^\triup(x)$  and $f^{\trido_1}(x) \subseteq f^{\trido_1}(1/2)$ for $0 \leq x \leq 1$.

DM1 can now choose an optimistic or a pessimistic approach. This decision depends on how much DM1 nows about the preference of DM2. If both have the same preference, then $\bar x = 1/2$ would be chosen by DM1, and DM2 would pick $(0,0)$ since $\{(0,0)\} = \MIN f^{\triup_1}(1/2)$ (the minimal points with respect to $\leq_{\R^2_+}$). If they have opposite preferences, the choice of $\bar x = 1/2$ would result in the outcome $\{(1,1)\} = \MAX f^{\triup_1}(1/2)$, i.e., a very bad one for DM1. Thus, in this case DM1 would pick either $x = 0$ or $x=1$ resulting in $(1,0)$ or $(0,1)$, respectively.

If DM2 has a different preference, the decision for DM1 is even more difficult. For example, if DM2 has the preference $R_2$ which is the vector order generated by the convex cone $C_2 = \{z \in \R^2 \mid z_1 \geq 0, z_2 \leq 0\}$ (positive $z_2$-values would be losses for DM2), then $f(x)$ only includes minimal points with respect to $R_2$ for all $x \in [0,1]$. This means DM1 has no clue which value in $f(x)$ DM2 is going to choose. The shared preference $R_1 \cap R_2$ is the vector order generated by the convex cone $C : = \{s (1,0) \mid s \geq 0\}$. Let $f^\triup$, $f^\trido$ denote the extensions of $f$ the extension of $f$ with respect to $R_1 \cap R_2$. The infimum of $f^\triup$ is $Q + C$  and the infimum of $f^\trido$ is $\vn$. Moreover, $f^\triup(1/2) \supseteq f^\triup(x)$ for all $x \in [1/2, 1]$ while the values $f^\triup(x)$ are not comparable with each other for $x \in [0, 1/2]$ and are all minimizers for $f^\triup$. Finally, the infimum of $f^\triup$ is attained in $[0, 1/2]$. Thus, if the decision makers collaborate, DM1 should choose $x \in [0, 1/2]$ where any outcome in $f^\triup(x)$ is possible since every point in $f^\triup(x)$ is minimal with respect to $\leq_C$. The situation does not change if DM2 decides with respect to $C_2$.

The situation changes if DM2 has a preference like the one generated by the cone $C'_2 = \{z \in \R^2 \mid z_2 \leq 0, z_2 \leq \frac{1}{2} z_1\}$. In this case, the minimal point of $f^{\triup_1}(x)$ with respect to $R'_2$ is $(1-2x, 1)$ for $x \in [0, 1/2]$ which might (should) provoke DM1 to choose the more pessimistic approach and minimize $f^{\trido_1}$ instead of $f^{\triup_1}$. 
\end{example}

\begin{example}
\label{ExDM2Punishes} Let $f \colon \R \to \P(\R_2)$ be given by
\[
f(x) =
\left\{
\begin{array}{ccc}
\cb{s(1, x) + (1-s)(0, x) \mid 0 \leq s \leq 1} & \colon & 0 \leq x \leq 1 \\
\vn & \colon & \text{otherwise}
\end{array}
\right.
\]
The relation $R_1$ is the vector order generated by the convex cone $C_1 = \{z \in \R^2 \mid z_1 \geq 0, z_2 \geq -z_1\}$ and $R_2$ is the vector preorder generated by $H^+(1,1) = \{z \in \R^2 \mid z_1 + z_2 \geq 0\}$. It is immediate that $f^{\triup_1}(0) \supseteq f^{\triup_1}(x)$ and $f^{\trido_1}(0) \subseteq f^{\trido_1}(x)$ for all $x \in X =[0,1]$. Therefore, $\{\bar x = 0\}$ produces the infimum both of $f^{\triup_1}$ and $f^{\trido_1}$ and is the best choice of DM1.

The final outcome is $(0,0)$ if DM2 decides according to $R_2$ or $R_1 \cap R_2 = R_1$ since $(0,0)$ is the only minimal point in $f(0)$ with respect to both $R_1$ and $R_2$. However, if DM2 does not want to get along with DM1 and decides according to $\hat R_ 2$, then all points in $f(0)$ are candidates for the final outcome since they are all minimal with respect to the vector order generated by $\hat C_2 = \{z \in \R^2 \mid z_1 \leq 0, z_2 \geq -z_1\}$.
\end{example}

As the previous examples show, DM1 has to decide which of $f^{\triup_1}$, $f^{\triup_2}$ or even $f^{\triup}$ should be minimized while DM2 faces a similar decision. In principle, the decision making procedures from Section \ref{SecSetOptDM} work where DM1 minimizes the chosen set-valued function while DM2 looks for minimizers of the function value selected by DM1 with respect to the chosen preference.

Note that $R \subseteq R_1, R_2$ and $\hat R_ i \subseteq R_i$, $i=1,2$. This means that there are potentially less comparable and more minimal 
elements for the smaller relations. Therefore, relationships between minimal elements and sets generating infima for subrelations are studied.

Note further that the two relations $\hat R_i$, $i = 1,2$, are not transitive in general, i.e., in many cases they are not preorders. As an example, one might consider the two vector orders on $\R^2$ generated by the two convex cones $\R^2_+$ and $H^+(1,1) = \cb{z \in \R^2 \mid z_1 + z_2 \geq 0}$. The set $(H^+(1,1) \bs \R^2_+) \cup\{0\}$ is a non-convex cone, thus the order relation generated by it is not transitive. 

The following assumption is made to remain in the framework of the present paper.

\smallskip
{\bf A.} Both $\hat R_1$ and $\hat R_2$ are transitive, i.e., preorders, whenever they are used.
\smallskip

A few details are added which show how the closure operators, the corresponding lattices and also minimal elements behave if sub-relations are considered. In the following, let $R_1 \subseteq R_2 \subseteq Z \times Z$ be two preorders on $Z$. The corresponding closure operators are denoted by $A^{\triup_i}$, $A^{\trido_i}$ and the lattices of closed sets by $\P^{\triup}(Z, \preceq_i)$, $\P^{\trido}(Z, \preceq_i)$, $i = 1,2$.

\begin{proposition}
\label{PropTwoClosures}
(1) One has $A^{\triup_1} \subseteq A^{\triup_2}$ for all $A \in \P(Z)$.

(2) One has $\of{A^{\triup_1}}^{\triup_2} = A^{\triup_2}$ for all $A \in \P(Z)$. 

(3) One has $\P^\triup(Z, \preceq_2) \subseteq \P^\triup(Z, \preceq_1)$ and $A = A^{\triup_1} = A^{\triup_2}$ for all $A \in \P^\triup(Z, \preceq_2)$.
\end{proposition}

{\sc Proof.} (1) If $z \in A^{\triup_1}$, then there is $a \in A$ with $a \preceq_1 z$. Since $R_1 \subseteq R_2$, this implies $a \preceq_2 z$, hence $z \in A^{\triup_2}$.

(2) On the one hand, one clearly has $A^{\triup_2} \subseteq \of{A^{\triup_1}}^{\triup_2}$ since $A \subseteq A^{\triup_1}$. On the other hand, (1) implies $\of{A^{\triup_1}}^{\triup_2} \subseteq \of{A^{\triup_2}}^{\triup_2} = A^{\triup_2}$ by idempotency of the closure operators.

(3) Assume $A \in \P^\triup(Z, \preceq_2)$, i.e., $A = A^{\triup_2}$. This implies $A \subseteq A^{\triup_1} \subseteq  A^{\triup_2} = A$ by (1) and the assumption, hence $A = A^{\triup_1}= A^{\triup_2}$. \pend

\begin{corollary}
\label{CorTwoClosureInfAttain}
(1) $\inf_{x \in X} f^{\triup_1}(x) \subseteq \inf_{x \in X} f^{\triup_2}(x)$. 

(2) If the infimum of $f^{\triup_1}$ is attained in $M \subseteq X$, then the infimum of $f^{\triup_2}$ is also attained in $M \subseteq X$.
\end{corollary}

{\sc Proof.} (1) This is immediate from (1) of Proposition \ref{PropTwoClosures} and the fact that the infimum is a union.

(2) Infimum attainment for $f^{\triup_1}$ means
\begin{equation}
\label{EqDelta1Inf}
\inf_{x \in X} f^{\triup_1}(x) = \bigcup_{x \in X} f^{\triup_1}(x) = \bigcup_{x \in M} f^{\triup_1}(x) = \inf_{x \in M} f^{\triup_1}(x).
\end{equation}
Proposition \ref{PropTriupUnion} yields $\bigcup_{x \in X} f^{\triup_1}(x) = \of{\bigcup_{x \in X} f(x)}^{\triup_1}$ and the same equation with $X$ replaced by $M$. Thus, by taking the $\triup_2$-closure, applying (2) of Proposition \ref{PropTwoClosures} and Proposition \ref{PropTriupUnion} again one gets
\[
\sqb{\of{\bigcup_{x \in X} f(x)}^{\triup_1}}^{\triup_2} = \of{\bigcup_{x \in X} f(x)}^{\triup_2} = \bigcup_{x \in X} f^{\triup_2}(x) = \inf_{x \in X} f^{\triup_2}(x)
\]
as well as the same equations with $X$ replaced by $M$. Using these identities, one gets $\inf_{x \in X} f^{\triup_2}(x) = \inf_{x \in M} f^{\triup_2}(x)$ by applying the $\triup_2$-closure to \eqref{EqDelta1Inf}. \pend

A parallel result for $f^{\trido_1}$, $f^{\trido_2}$ does not seem possible since the infimum with respect to $\subseteq$ is an intersection which does not commute with the closure operators (see comment before Lemma \ref{LemLowerSupMax}). The following examples illustrate the situation.

\begin{example}
\label{ExTwoRelationInfimum}
(1) Consider the function $f \colon \R \to \P(\R^2)$ given by
\[
f(x) =
\left\{
\begin{array}{ccc}
\co\cb{(0,0)^\top, (x, -x)^\top} & \colon & 0 \leq x \leq 1 \\
\vn & \colon & \text{otherwise}
\end{array}
\right.
\]
Let $R_1$ be the vector order generated by the convex cone $\R^2_+$ and $R_2$ the vector preoder generated by the half-space $H^+(1,1)$. Then $M = \{1\}$ produces the infimum for $f^{\triup_1}$ and also for  $f^{\triup_2}$. However, every $x \in [0,1]$ produces the infimum of $f^{\triup_2}$ since $f^{\triup_2}(x) = H^+(1,1)$ for all $x \in [0,1]$.

(2) Consider the function $f \colon \R \to \P(\R^2)$ given by
\[
f(x) =
\left\{
\begin{array}{ccc}
\cb{(x, -2x+3)^\top} & \colon & 0 \leq x \leq 1 \\
\cb{(x, -\frac{1}{2}x+\frac{3}{2})^\top} & \colon & 1 \leq x \leq 3 \\
\vn & \colon & \text{otherwise}
\end{array}
\right.
\]
and the same relations as in (1). Then $M = [0,3]$ is the smallest set producing the infimum of $f^{\triup_1}$. On the other hand, $N = \{1\}$ already produces the infimum of $f^{\triup_2}$.
\end{example}

Note that $\inf_{x \in \R} f^{\triup_1}(x) \subsetneq \inf_{x \in \R} f^{\triup_2}(x)$ for (1), (2) of Example \ref{ExTwoRelationInfimum}, i.e., the two infima are different in both cases. 

The major consequence of Corollary \ref{CorTwoClosureInfAttain} is that if the infimum of the lattice extension of a function is attained in a set with respect to the relation $R_1 \cap R_2$, then it is also attained both with respect to $R_1$ and $R_2$. Similarly, if it is attained with respect to $\hat R_1$, then it is also attained with respect to $R_1$.

The situation is different for minimal points.

The decision makers could be interested in finding elements which cannot be improved with respect to both preferences simultaneously, i.e., with respect to $R_1 \cap R_2$. In the following, the minimal points of $A \subseteq Z$ with respect to $R_1$, $R_2$ and $R_1 \cap R_2$ are denoted by $\MIN_1 A$, $\MIN_2 A$ and $\MIN_{12} A$, respectively. The next lemma is needed to ensure consistency of the sets of minimal points. It should be known but a reference could not be found by the author.

\begin{lemma}
\label{LemMinMonotonicity}
Let $R_1, R_2$ be two preorders on $Z$ with $R_1 \subseteq R_2$. If $S(R_2) \cap R_1 \subseteq S(R_1)$, then $\MIN_2 A \subseteq \MIN_1 A$ for $A \subseteq Z$. 
\end{lemma}

{\sc Proof.} Take a point $\bar z \in \MIN_2 A$. If $(z, \bar z) \in R_1$, then $(z, \bar z) \in R_2$ since $R_1 \subseteq R_2$. This implies $(z, \bar z) \in S(R_2)$ since $\bar z \in \MIN_2 A$. Thus $(z, \bar z) \in S(R_1)$ by assumption, hence $\bar z \in \MIN_1 A$. \pend

Note that the condition $S(R_2) \cap R_1 \subseteq S(R_1)$ is automatically satisfied if $R_2$ is a partial order since in this case its symmetric part is just the diagonal relation.

Applying the lemma to $R_1 \cap R_2$ and $R_1$ (and $R_2$) one obtains the following.

\begin{corollary}
\label{CorMinIntersection}
If $S(R_2) \cap R_1 \subseteq S(R_1 \cap R_2)$ and $S(R_1) \cap R_2 \subseteq S(R_1\cap  R_2)$ is satisfied, then $\MIN_1 A \cap \MIN_2 A \subseteq \MIN_{12} A$.
\end{corollary}

{\sc Proof.} Replacing $R_1$ in Lemma \ref{LemMinMonotonicity} by $R_1 \cap R_2$ and $R_2$ by $R_1$ one gets $\MIN_1 A \subseteq \MIN_{12} A$ whenever $S(R_1) \cap (R_1 \cap R_2) \subseteq S(R_1 \cap R_2)$. Clearly, $S(R_1) \cap (R_1 \cap R_2) = S(R_1) \cap R_2$. The parallel condition is obtained for $R_2$. \pend

Thus, the decision makers can be sure that elements in $\MIN_1 A \cap \MIN_2 A$ cannot be improved by with respect to either preference. However, this set can be empty as shown in the following example.

\begin{example}
(1) Let $A = \{z \in \R^2 \mid z_2 \geq -2z_1 + 2, 2z_2 \geq -z_1 + 2, z_1 \geq 0, z_2 \geq 0\}$ and $R_1$ the partial order generated by the convex cone $\R^2_+$, $R_2$ the total preorder generated by the cone $\{z \in \R^2 \mid z_2 \geq -z_1\}$ which is a half-space. Then $\MIN_1 A$ is the union of the two line segments between $(2,0)$, $(1,1)$ and $(1,1)$, $(0,2)$, respectively, while $\MIN_2 A$ only includes the point $(1,1)$. Clearly, $\MIN_2 A \subseteq \MIN_1 A$ and $S(R_2) \cap R_1 \subseteq S(R_1)$ is satisfied.

(2) Let $A$ and $R_1$ the same as in (1) and $R_2$ be the complete preorder generated by the half-space $\{z \in \R^2 \mid z_1 \geq 0\}$. Then $\MIN_2 = \{(0, z_2) \mid z_2 \geq 2\}$ which is not a subset of $\MIN_1 A$. Of course, $S(R_2) \cap R_1 \subseteq S(R_1)$ is not satisfied but one has $\MIN_1 A \cap \MIN_2 A = \{(0,2)\} \neq \vn$.

(3) Let $A \subseteq \R^2$ now be the square with vertices $(0,0)$, $(0,1)$, $(1,0)$, $(1,1)$ and $R_1$ the complete preorder from (2) while $R_2 = (R_1)^{-1}$. Then $\MIN_1 A$ is the line segment between $(0,0)$, $(0,1)$ while $\MIN_2 A$ is the one between $(1,0)$, $(1,1)$, thus $\MIN_1 A \cap \MIN_2 A = \vn$. Moreover, $R_1 \cap R_2 = S(R_1) = S(R_2)$ which shows that the two conditions of Corollary \ref{CorMinIntersection} are satisfied. One has $\MIN_{12} A = A$. Of course, this is an example for which $R_2 = (R_1)^{-1}$ (conflicting preorders) but $R_1 \cap R_2 \neq \vn$.
\end{example}

By the way of conclusion, the decision making procedure in the case of only partially aligned preferences highly depends on the knowledge (or the belief) of the decision makers about the behaviour of the respective opponent. Results like Corollary \ref{CorTwoClosureInfAttain} for infimum attainment and Lemma \ref{LemMinMonotonicity} for minimal elements show a behaviour in different directions with respect to relations $R_1$, $R_2$ with $R_1 \subseteq R_2$. This is in accordance with the features of solution sets for set optimization problems: a subset $M$ of $X$ should be big enough to produce the infimum and small enough such that $\{f(x) \mid x \in M\}$ only includes minimal function values: see Definitions \ref{DefSolConceptColl}, \ref{DefSolConceptConf}.

\section{An application to games with vector payoffs}
\label{SecGames}

Let $Y$ be another non-empty set. A function $F_1 \colon X \times Y \to Z$ is called a payoff (loss) function for the game where player P1 has to select $x \in X$ and player P2 $y \in Y$ resulting in the loss $F_1(x, y)$ for P1. P1 has the preference $\preceq_1$ and wants to minimize losses. If P1 chooses $x \in X$, then anything in $f_1(x) = \{F_1(x, y) \mid y \in Y\}$ can be the loss. The worst case for P1 is a choice $y \in Y$ of P2 such that $F_1(x, y) \in \MAX f_1(x)$.

Thus, a worst case analysis of such a two-person game would be to consider it as a conflicting decision making procedure: P1 as DM1 wants to minimize $f_1$ and P2 is DM2 with inverse preference and a corresponding loss function $F_2$.

The same situation occurs if $F_1(x, y)$ is considered as uncertain outcome which varies with $y \in Y$; DM1 wants to make a robust decision, i.e., one that minimizes the worst possible loss.

Defining $f_1^\trido(x) = \{z \in Z \mid \exists y \in Y \colon z \preceq_1 F_1(x, y)\}$ one gets the complete lattice-valued version of the problem for P1:
\[
\text{minimize} \quad f_1^\trido(x) \quad \text{subject to} \quad x \in X.
\]
If player P2 has the loss function $F_2 \colon X \times Y \to Z$ one defines $f_2$, $f_2^\trido$ in a parallel fashion as functions of $y \in Y$ arriving at the problem
\[
\text{minimize} \quad f_2^\trido(y) \quad \text{subject to} \quad y \in Y.
\]

\begin{definition}
\label{DefWCEquilibrium}
A pair $(\bar x, \bar y) \in X \times Y$ is called a worst-case equilibrium for the game $(X, Y, F_1, F_2)$ if $f_1^\trido(\bar x)$ is minimal in $\cb{f^\trido_1(x) \mid x \in X}$ and $f_2^\trido(\bar y)$ is minimal in $\cb{f^\trido_2(y) \mid y \in Y}$, both with respect to $\subseteq$. 
\end{definition}

Worst-case equilibria are interchangeable by definition: if $(\bar x, \bar y)$ and $(\bar x', \bar y')$ are worst-case equilibria, so are $(\bar x, \bar y')$ and $(\bar x', \bar y)$. 

On the other hand, if P1 knows or assumes that P2 plays $\bar y \in Y$, a best response would be $\bar x \in X$ with $F_1(\bar x, \bar y) \in \MIN\{F_1(x, \bar y) \mid x \in X\}$. The function $g_1(y) = \{F_1(x, y) \mid x \in X\}$ is the set of possible outcomes as response of P1 to P2's choice $y \in Y$. Its $\P^\triup(Z, \preceq_1)$-valued extension is defined by $g^\triup_1(y) = \{z \in Z \mid \exists x \in X \colon F_1(x, y) \preceq_1 z\}$, so that $g_1(y)$ and $g^\triup_1(y)$ have the same set of minimal elements. Likewise, $g_2(x)$ and $g^\triup_2(x)$ are defined.

\begin{definition}
\label{DefShapleyEquilibrium}
A pair $(\bar x, \bar y) \in X \times Y$ is called a Shapley equilibrium for the game $(X, Y, F_1, F_2)$ if 
\[
F_1(\bar x, \bar y) \in \MIN g^\triup_1(\bar y) \quad \text{and} \quad F_2(\bar x, \bar y) \in \MIN g^\triup_2(\bar x).
\]
\end{definition}

Shapley equilibria are not interchangeable in general as already an example for a zero-sum matrix game with vector payoff shows: see \cite{HamelLoehne18MMOR}, the references therein and also \cite{Park19JTE}. Moreover, vector games can have worst-case equilibria whose payoff is not worse than the payoff of a Shapley equilibrium. This phenomenon already occurs for bi-matrix games, see \cite{AumannMaschler72MS, Park19JTE, Pruzhansky11IJGT}. Thus, it is not clear what to recommend players of such games. Finally, there can be any number of Shapley equilibria which also makes the choice difficult. 

\begin{remark}
\label{RemSetOptGame}
On the one hand, P1 wants to minimize the function $f^\trido_1 \colon X \to \P^\trido(Z, \preceq_1)$ which has the "wrong" image lattice. One the other hand, P1 wants to have minimal elements of values of the function $g^\triup_1 \colon Y \to \P^\triup(Z, \preceq_1)$ which has the "right" image lattice for minimization but P1 cannot influence which of the values of $g^\triup_1$ is realized.

Therefore, the decision procedure cannot be a simple application of DPConf1.
\end{remark}

Focusing on Shapley equilibria, the decision procedure for P1 could just be to select $\bar x \in X$ such that there is $\bar y \in Y$ forming the Shapley equilibrium $(\bar x, \bar y)$. With the issues mentioned above in view one can design a decision making procedure taking worst case equilibria into account.

\medskip
\noindent
{\bf Decision Procedure for Games (DPGame)} (for P1, P2 is parallel)
\begin{description}
\item[{\bf 1}] P1 determines a set $Z_1 \subseteq Z$ such that for each $z \in Z_1$ there is $y \in Y$ with $f^\trido_2(y)$ is minimal and $z \in  \MIN g^\triup_1(y)$. 
\item[{\bf 2}] P1 determines $\bar x \in X$ such that $f^\trido_1(\bar x)$ is minimal and $Z_1 \subseteq f^\trido_1(\bar x)$. 
\end{description}

Of course, there still is plenty of "free will," i.e., there might be several, even many elements in $Z_1$ and many such $\bar x$, and there is no guarantee the opponent follows suit. Finally, similar decision procedures can be designed starting with sets which are merely acceptable.

\begin{remark}
A reviewer of the first version of this paper pointed out that leader-follower games (see, for example, \cite{VonStengelZamir10GEB}) might even better fit into the proposed framework than the game situation with simultaneous decisions. Clearly, the two-person game above can be considered as an as-if two-stage decision procedure since each player anticipates the moves of the opponent and reacts in advance to them. The two solution concepts in Definition \ref{DefWCEquilibrium} and \ref{DefShapleyEquilibrium} reflect the fact that both players can act as if they were first or second stage decision makers--which makes games difficult. The discussion of leader-follower games, especially with incomplete preferences for the players, goes beyond the scope of this paper and is worth a separate study.
\end{remark}

\section{Conclusions}
\label{SecConclusions}

In this paper, a methodology is developed which turns the 'preference for flexibility' approach into two-stage decision making procedures via set optimization methods. There are already applications of this type such as discussed in \cite{HamelLoehne24ArX}. The case of completely aligned preferences and of inverse preferences for the two stages can be dealt with.

The case of partially aligned preferences seems to deserve more attention since it appears to be realistic that in conflict situations the opponents at least agree on comparisons in certain directions and only disagree on others; a fundamental question might be how to steer the negotiations in such a case towards collaboration. In particular, games with partially aligned preferences seem to be an unexplored area. A comparison of different approaches for this case may try to quantify the difference of the preferences, i.e., how ``big" $\hat R_ 1$, $\hat R_ 2$ in Section \ref{SecPartiallyAligned} are. Similarly, one can study the impact of a more collaborative and a more contradictive approach to the actual outcome for the decision makers. The situation resembles the prisoners' dilemma now at the level of choosing a preference.

It seems that several applications fit naturally into the proposed framework, among them leader-follower games as pointed out by a reviewer and also deterministic bilevel optimization problems where the upper level decision has to be made in advance (or without the knowledge) of the lower level decision.

Finally, the case can be considered when the first stage decision maker has no or only partial knowledge of the preference of second stage one.

\bigskip\noindent
{\bf Acknowledgement.} The author is extremely grateful to reviewers who provided a number of constructive comments which not only led to an improvement of the paper, but also hinted at potential extensions. Many thanks go to Andreas L\"ohne with whom the author had several inspiring discussions on topics related to this article.



\end{document}